\documentclass[runningheads]{llncs}
\usepackage[T2A]{fontenc}
\usepackage[utf8]{inputenc}
\usepackage{url}
\usepackage{multirow}
\usepackage{amssymb}
\usepackage{amsxtra}
\usepackage{makeidx}
\usepackage{amsfonts}
\usepackage{algorithm}
\usepackage{algorithmic}
\usepackage{float}
\usepackage{graphicx}
\usepackage[font=small,labelfont=bf]{caption}
\usepackage{comment}
\usepackage{color}
\usepackage{hyperref}
\usepackage{nicefrac}
\usepackage{subcaption}

\begin{document}
\title{On Modification of an Adaptive Stochastic Mirror Descent Algorithm for Convex Optimization Problems with Functional Constraints\thanks{This paper accepted to the print as a chapter in the forthcoming book: \textit{Communications in Mathematical Computations and Applications, IACMC2019, Springer.}}}
%
%\titlerunning{Proc. International Arab Conference on Mathematics and Computations 2019(IACMC2019), 1-7}
\author{Mohammad S. Alkousa \orcidID{0000-0001-5470-0182} }
\authorrunning{M.~S.~Alkousa}
\institute{Moscow Institute of Physics and Technology, Moscow, Russia\\
	\email{ mohammad.alkousa@phystech.edu}	}

\maketitle
\begin{abstract}
This paper is devoted to a new modification of a recently proposed adaptive stochastic mirror descent algorithm for constrained convex optimization problems in the case of several convex functional constraints. Algorithms, standard and its proposed modification, are considered for the  type of problems with non-smooth Lipschitz-continuous convex objective function and convex functional constraints. Both algorithms, with an accuracy $\varepsilon$ of the approximate solution to the problem, are optimal in the terms of lower bounds of estimates and have the complexity $O\left( \varepsilon^{-2} \right)$. In both algorithms, the precise first-order information, which connected with (sub)gradient of the objective function and functional constraints, is replaced with its unbiased stochastic estimates. This means that in each iteration, we can still use the value of the objective function and functional constraints at the research point, but instead of their (sub)gradient, we calculate their stochastic (sub)gradient. Due to the consideration of not all functional constraints on non-productive steps, the proposed modification allows saving the running time of the algorithm. Estimates for the rate of convergence of the proposed modified algorithm is obtained. The results of numerical experiments demonstrating the advantages and the efficient of the proposed modification for some examples are also given.

\end{abstract}

\section{Introduction}
Large scale non-smooth convex optimization is a common problem for a range of computational areas including statistics, computer vision, general inverse problems, machine learning, data science and in many applications arising in applied sciences and engineering. Since what matters most in practice is the overall computational time to solve the problem, first-order methods with computationally low-cost iterations become a viable choice for large scale optimization problems.

Generally, first-order methods have simple structures with a low memory requirement. Thanks to these features, they have received much attention during the last decade. There are a lot of first-order methods for solving the optimization problems in the case of non-smooth objective function. Some examples of these methods, to name but a few, are: subgradient methods \cite{Nesterov_2018_book,polyak_book,shor_book}, subgradient projection methods \cite{Nesterov_2018_book,polyak_book,shor_book}, OSGA \cite{Neumaier}, bundle-level method \cite{Nesterov_2018_book}, Lagrange multipliers method \cite{bib_Boyd} and many others.

There is a long history of studies on continuous optimization with functional constraints. The recent works on first-order methods for convex optimization with convex functional constraints include \cite{bib_Adaptive,paper:almost_surely_constrained_2019,paper:Lin_2018,article:titov_optima_2019,paper:Wei_2018,paper:Xu_2019} for deterministic constraints and \cite{article:mohammad_kim_2019,paper:Basu_2019,paper:Lan_2016,paper:Xu_primal_dual_2019} for stochastic constraints.  However, the parallel development for problems with non-convex objective functions and also with non-convex constraints, especially for theoretically provable algorithms, remains limited, see \cite{paper:nonconvex_obj_and_cons_2019} and references therein.

The mirror descent algorithm which originated in \cite{nemirovsky1983problem,nemirovskii1979efficient} and  was later analyzed in \cite{beck2003mirror}, is considered as the non-Euclidean extension of subgradient methods. The standard subgradient methods employ the Euclidean distance function with a suitable step-size in the projection step. Mirror descent extends the standard projected subgradient methods by employing a nonlinear distance function with an optimal step-size in the nonlinear projection step \cite{mirror_weight}.  Mirror descent method not only generalizes the standard gradient descent method, but also achieves a better convergence rate \cite{article:doan_2019}. In addition, Mirror descent method is applicable to optimization problems in Banach spaces where gradient descent is not \cite{article:doan_2019}. An extension of the mirror descent method for constrained problems was proposed in \cite{beck2010comirror,nemirovsky1983problem}. 

Usually, the step-size and stopping rule for mirror descent algorithms require to know the Lipschitz constant of the objective function and constraint, if any. Adaptive step-sizes, which do not require this information, are considered for unconstrained problems in \cite{bib_Nemirovski}, and for constrained problems in \cite{beck2010comirror}. Some optimal mirror descent algorithms, for convex optimization problems with non-smooth convex functional constraint and both adaptive step-sizes and stopping rules, are proposed in \cite{bib_Adaptive}.  Also, there were considered some modifications of these algorithms for the case of problems with many functional constraints in \cite{bib_Stonyakin}.

If we focus on the problems of minimization of an objective function consisting of a large number of component functionals, such as $f(x) = \sum\limits_{j=1}^{N} f_j(x)$ where $f_j : \mathbb{R}^n \rightarrow  \mathbb{R}, j = \overline{1,N}$ are convex, then in each iteration of any iterative minimization procedure  computing a single (sub)gradient $\nabla f(x) = \sum\limits_{j=1}^{N} \nabla f_j(x)$ becomes very expensive. Therefore there is an incentive to calculate the stochastic (sub)gradient $\nabla f(x, \zeta)$ where $\zeta$ is a random variable taking its values in $\{1,\ldots,N\}$. This mean that $\nabla f(x,\zeta) = \nabla f_i(x)$, were $i$ is chosen randomly in each iteration from the set $\{1,\ldots,N\}$, or instead, one can  employ randomly chosen a mini-bach approach in which a small subset $S \subset \{1,\ldots,N\}$ is chosen randomly, then  $\nabla f(x, \zeta) = \sum\limits_{i \in S} \nabla f_i(x)$. This randomly calculating of the (sub)gradient is known as stochastic (sub)gradient.

In the stochastic version of an optimization method, the exact first-order information is replaced with its unbiased stochastic estimates, where the exact first-order information is unavailable. This permits accelerating the solution process, with the earning from randomization growing progressively with problem’s sizes. A different approach to solving stochastic optimization problems is called stochastic approximation (SA), which was initially proposed in a seminal paper by Robbins and Monro in 1951 \cite{Robbins_Monro_1951}. An important improvement of this algorithm was developed by Polyak and Juditsky \cite{Polyak_Judsky,Polyak_stoc}. More recently, Nemirovski et al. \cite{Nemirovski_mirror_stoc} presented a modified stochastic approximation method and demonstrated its superior numerical performance for solving a general class of non-smooth convex problems.

This paper is devoted to a new modification of an adaptive stochastic mirror descent algorithm (see Algorithm 4 in \cite{bib_Adaptive}. This algorithms is listed as Algorithm \ref{algorithm1}, below), which is proposed to solve the stochastic setup (randomized version) of the convex minimization problems in the case of several convex functional constraints. This means that we can still use the value of the objective function and functional constraints at the research point, but instead of their (sub)gradient, we use their stochastic (sub)gradient. Namely, that we consider the first-order unbiased oracle that produces  stochastic (sub)gradients of the objective function and functional constraints, see for example \cite{Duchi_2016,Shpiro_2014}. We consider the arbitrary proximal structure and the type of problems with non-smooth Lipschitz-continuous objective function.  Furthermore, it has been proved a theorem to estimate the rate of convergence of the proposed modification, from this theorem we can see that the modified algorithm achieves the optimal complexity of the order $O \left(\varepsilon^{-2}\right)$ for the class of problems under consideration (see \cite{nemirovsky1983problem}). 

The rest of the paper is organized as follows. In Section \ref{problem_basics} we give some basic notation, summarize the problem statement and standard mirror descent basics. In Section \ref{Adaptive_SMD}  we display the adaptive stochastic mirror descent algorithm (Algorithm 4 in \cite{bib_Adaptive}). Section \ref{section_of_mod} is devoted to the proposed modified algorithm and proving a theorem about the rate of convergence of this algorithm and its optimal complexity estimate. In the last section, we consider some numerical experiments that allow us to compare the work of standard algorithm and its proposed modification for certain examples.

\section{Problem Statement and Standard Mirror Descent Basics}\label{problem_basics}
Let $\mathbb{V}$ be a finite-dimensional vector space, endowed with the norm $\|\cdot\|$, and $\mathbb{V}^*$ is the conjugate space of $\mathbb{V}$ with the following norm
$$
	\|h\|_*=\max\limits_x\{\langle h,x\rangle,\|x\|\leq1\},
$$
where $\langle h,x\rangle$ is the value of the continuous linear functional $h$ at $x \in \mathbb{V}$.

Let $Q\subset \mathbb{V}$ be a closed convex set, $f$ and $g_j : Q \rightarrow \mathbb{R} \, (j = \overline{1, m})$  convex subdifferentiable functionals. We assume that $f$ and $g_j (j = \overline{1, m})$ are Lipschitz-continuous, i.e. there exist $M_f > 0$ and  $M_g>0$, such that
\begin{equation}
  |f(x)-f(y)|\leq M_f\|x-y\|\quad \forall \; x,y\in Q,
\end{equation}
\begin{equation}
|g_j(x)-g_j(y)|\leq M_g\|x-y\|\quad \forall \; x,y\in Q, \; j = \overline{1, m}.
\end{equation}

It is clear that instead of a set of functionals  $\{g_j(\cdot)\}_{j=1}^{m}$ we can see one functional $g: Q \rightarrow \mathbb{R}$, such that
$$
	g(x) = \max\limits_{j = \overline{1, m}} \{g_j(x)\}, \quad |g(x)-g(y)|\leq M_g\|x-y\|\; \quad \forall \; x,y\in Q. 
$$
It means that at every point $x \in Q$ there is a subgradient $\nabla g(x)$, and $\|\nabla g(x)\|_* \leq M_g$. Recall that for a differentiable functional $g$, the subgradient $\nabla g(x)$ coincides with the usual gradient.

In this paper, we consider the stochastic setup of the following convex constrained optimization problem
\begin{equation}\label{problem}
f(x) \rightarrow \min\limits_{x\in Q, \; g(x) \leq 0}.
\end{equation}

For the stochastic setup of the problem \eqref{problem}, we introduce the following assumptions (see \cite{bayandina_stoc,bib_Adaptive}). Given a point $x \in Q$, we can calculate the stochastic (sub)gradients $\nabla f(x,\xi)$ and $\nabla g(x,\zeta)$, where $\xi$ and  $\zeta$ are random vectors. These stochastic (sub)gradients satisfy
\begin{equation}\label{bound_stoc_f_g_1}
\mathbb{E}[\nabla f(x,\xi)] = \nabla f(x) \in \partial f(x) \quad \text{and} \quad  \mathbb{E}[\nabla g(x,\zeta)] = \nabla g(x) \in \partial g(x),
\end{equation}
where $\mathbb{E}$ denote to the expectation, and
\begin{equation}\label{bound_stoc_f_g_2}
\|\nabla f(x,\xi)\|_* \leq M_f \quad  \text{and} \quad \|\nabla g(x,\zeta)\|_* \leq M_g, \;\;\; a.s. \; \text{in} \; \xi, \zeta.
\end{equation}

To motivate these assumptions, let $ S_n(1) = \left\{ x\in \mathbb{R}_{+}^n \;|\; \displaystyle \sum_{i=1}^{n}x_i = 1\right\} $ be a standard unit simplex in $\mathbb{R}^n$, we consider the following optimization problem
\begin{equation*}
\begin{cases}
f(x) = \frac{1}{2} \langle Ax,x \rangle  \rightarrow \min\limits_{x \in S_n(1)},\\
s.t. \ \  g(x)= \max\limits_{i = \overline{1, m}} \{\langle c_i, x \rangle\}  \leq 0,
\end{cases}
\end{equation*}
where $A$ is a given $n \times n$ matrix and  $c_i \; (i=\overline{1, m})$ are given vectors in $\mathbb{R}^n.$ (See \cite{bib_Adaptive})\\

The exact computation of the gradient $\nabla f(x) = Ax$ takes $O(n^2)$ arithmetic operations, which is expensive, when $n$ is very large, for the huge-scale optimization problems. In this setting, it is natural to use the randomization to construct a stochastic approximation for  $\nabla f(x)$. Let $\xi$ be a random variable its values $1, \ldots, n$ with probabilities $x_1, \ldots, x_n$ respectively. Let $A^{\langle i \rangle}$ denote the $i$-th column
of the matrix $A$. Since $x \in S_n(1)$,
$$
	\begin{aligned}
	\mathbb{E}[A^{\langle \xi \rangle}] &= A^{\langle 1 \rangle} \underbrace{\mathbb{P}(\xi =1)}_{x_1} + \cdots + A^{\langle n \rangle} \underbrace{\mathbb{P}(\xi =n)}_{x_n}\\
	& = A^{\langle 1 \rangle} x_1 + \cdots + A^{\langle n \rangle} x_n = Ax,
	\end{aligned}
$$
where $\mathbb{P}$ denote to the probability of an event.

Thus, we can use $A^{\langle \xi \rangle}$  as a stochastic gradient of $f$ (i.e. $\nabla f(x,\xi ) = A^{\langle \xi \rangle}$), which can be calculated in $O(n)$
arithmetic operations.

Let $d : Q \rightarrow \mathbb{R}$ be a distance generating function, which is continuously differentiable and $1$-strongly convex with respect to the norm $\|\cdot\|$, i.e.
$$
	 \hspace{0.2cm} d(y) \geq d(x) +\langle \nabla d(x), y-x \rangle + \frac{1}{2} \| y-x \|^2 \quad \forall\;  x, y \in Q,
$$
and assume that $\min\limits_{x\in Q} d(x) = d(0).$ Suppose, we have a constant $\Theta_0 >0$ such that $d(x_{*}) \leq \Theta_0^2,$ where $x_*$ is a solution to the problem \eqref{problem}.

Note that if there is a set of optimal points for \eqref{problem} $X_* \subset Q$,  we may assume that
$$
	\min\limits_{x_* \in X_*} d(x_*) \leq \Theta_0^2.
$$

For all $x, y\in Q \subset \mathbb{V}$, we consider the corresponding Bregman divergence,  which was initially studied by Bregman \cite{Bregman} and later by many others (see \cite{Bauschke_et.al}),
$$
	V_x(y) = d(y) - d(x) - \langle \nabla d(x), y-x \rangle.
$$

In particular, in the standard proximal setup (i.e. Euclidean setup) we can choose $d(x) = \frac{1}{2} \|x\|_2^2$, then $V_x(y) = \frac{1}{2} \|x-y\|_2^2$. Another setups, for example entropy, $\ell_1/\ell_2$, simplex, spectahedron and many others, can be found in \cite{bib_Nemirovski}.

We also assume that the constant $\Theta_0>0$ is known, such that
\begin{equation}\label{cond_on_Bregman}
\sup\limits_{x,y \in Q} V_x(y) \leq \Theta_0^{2}.
\end{equation}

For all  $x\in Q $ and $p\in \mathbb{V}^*$, the proximal mapping operator (mirror descent step) is defined as
$$
	\mathrm{Mirr}_x (p) = \arg\min\limits_{u\in Q} \big\{ \langle p, u \rangle + V_x(u) \big\}.
$$

We make the simplicity assumption, which means that $\mathrm{Mirr}_x (p)$ is easily computable.

Let $x_*$ be a solution to \eqref{problem} and $\varepsilon > 0$ is given, we say that a (random) point $\hat{x} \in Q$ is an expected $\varepsilon$-solution to \eqref{problem} if
\begin{equation}\label{expected_sol}
\mathbb{E}[f(\hat{x})] - f(x_*) \leq \varepsilon \quad \text{and}\quad g(\hat{x}) \leq \varepsilon. 
\end{equation}

The following well-known lemma describes the main property of the proximal mapping operator (see  \cite{bib_Adaptive,bib_Nemirovski}).

\begin{lemma}\label{lem1}
	Let $f: Q \rightarrow \mathbb{R}$ be a convex subdifferentiable function over the convex set $Q$ and  $z=Mirr_{y}\left(h \nabla f(y, \xi)\right)$ for some $h>0$, $y, z \in Q$ and $\xi$ random vector. Then for each $x\in Q$ we have
	\begin{equation*}\label{eq7}
		h\left(f(y) - f(x) \right)\leq\frac{h^2}{2}\|\nabla f(y,\xi)\|_*^2 + V_y(x) - V_z(x) + h \left\langle \nabla f(y,\xi)-\nabla f(y), y-x \right\rangle.
	\end{equation*}
\end{lemma}

\section{Adaptive Stochastic Mirror Descent Algorithm}\label{Adaptive_SMD}

In \cite{bib_Adaptive} it was considered an adaptive method, for the convex optimization problem \eqref{problem} in the stochastic setup described above (see Algorithm \ref{algorithm1}). In this setting, the output of the algorithm is random, in the sense of \eqref{expected_sol}. The adaptivity of this method is in terms of step-size and stopping role, which is mean that we do not need to know the constants $M_f$ and $ M_g$ in advance. We assume that, on each iteration of the algorithm, independent realizations of the random variables $\xi$ and $\zeta$ are generated. In this section, we show this algorithm and the fundamental result of the estimate about the convergence rate of this algorithm.

As can be seen from the items of the Algorithm \ref{algorithm1}, the needed point (Ensure) is selected among the points $x^k$ for which $g(x^k) \leq \varepsilon$. Therefore, we will call step $k$ productive if $g(x^k) \leq \varepsilon$. If the reverse inequality $g(x^k) > \varepsilon$ holds then step $k$ will be called non-productive.

Let $I, J$ denote the set of indexes of productive  and non-productive steps produced by Algorithm \ref{algorithm1}, respectively. $N_I, N_J$ denote the number of productive and non-productive steps, respectively.

\begin{algorithm}[t]
	\caption{Adaptive Stochastic Mirror Descent Algorithm.}
	\label{algorithm1}
	\begin{algorithmic}[1]
		\REQUIRE  accuracy $\varepsilon,$ starting point $x^0, d(\cdot),  Q, \Theta_0$ such that \eqref{cond_on_Bregman} holds.
		\STATE $I=:\emptyset$
		\STATE $N\leftarrow0$
		\REPEAT
		\IF{$g(x^N) \leq \varepsilon$}
		\STATE $M_N := \left\|\nabla f(x^N,\xi^{N})\right\|_*,$
		\STATE $h_N = \Theta_0 \left(\sum\limits_{t = 0}^N M_t^2\right)^{-1/2},$
		\STATE $x^{N+1}:=Mirr_{x^N}\left(h_N\nabla f(x^N,\xi^{N})\right),$  "productive step"
		\STATE$N\rightarrow I$
		\ELSE
		\STATE $M_N := \left\|\nabla g(x^N,\zeta^N)\right\|_*,$
		\STATE $h_N = \Theta_0 \left(\sum\limits_{t = 0}^N M_t^2\right)^{-1/2},$
		\STATE $x^{N+1}:=Mirr_{x^N}\left(h_N\nabla g(x^N,\zeta^N)\right),$    "non-productive step"
		\ENDIF
		\STATE $N\leftarrow N+1$
		\UNTIL{$N  \geq \frac{2 \Theta_0}{\varepsilon}\left(\sum\limits_{t = 0}^{N-1} M_t^2\right)^{1/2}.$}
		\ENSURE $\bar{x}^N:=\frac{1}{N_I}\sum\limits_{k\in I}x^k.$
	\end{algorithmic}
\end{algorithm}

For the complexity estimate of Algorithm \ref{algorithm1}, the next result was obtained in  \cite{bayandina_stoc,bib_Adaptive}.
\begin{theorem}\label{theo_standard}
Let equalities \eqref{bound_stoc_f_g_1} and inequalities \eqref{bound_stoc_f_g_2} hold. Assume that a known constant $\Theta_0>0$ is such that inequality \eqref{cond_on_Bregman} holds. Then Algorithm \ref{algorithm1} stops after no more than\\
\begin{equation}\label{estimate_alg1}
N = \left\lceil \frac{4 \max\{M_f^2, M_g^2\} \Theta_0^2}{\varepsilon^2}\right\rceil
\end{equation}
iterations and $\bar{x}^N$ is an expected $\varepsilon$-solution to problem \eqref{problem} in the sense of \eqref{expected_sol}.
\end{theorem}

\section{The Modification of an Adaptive Stochastic Mirror Descent Algorithm }\label{section_of_mod}

In this section, we consider a modification of an Algorithm \ref{algorithm1}. The idea of this modification was considered in \cite{bib_Stonyakin} for some adaptive mirror descent algorithms to solve the deterministic setup of the convex optimization problems with Lipschitz-continuous functional constraints. This idea is summarized as: when we have a non-productive step $k$, i.e. $g(x^k)> \varepsilon$, then instead of calculating the subgradient of the functional constraint with max-type $g(x) = \max\limits_{i = \overline{1, m}} \{g_i(x)\}$, we calculate (sub)gradient of one functional $g_j$, for which we have $g_j(x^k)> \varepsilon$. The proposed modification allows saving the  running time of algorithm due to consideration of not all functional constraints on non-productive steps.

\begin{algorithm}[H]
	\caption{ The Modification of an Adaptive Stochastic Mirror Descent Algorithm.}
	\label{algorithm2}
	\begin{algorithmic}[1]
		\REQUIRE  accuracy $\varepsilon,$ starting point $x^0, d(\cdot),  Q, \Theta_0$ such that \eqref{cond_on_Bregman} holds.
		\STATE $I=:\emptyset$
		\STATE $N\leftarrow0$
		\REPEAT
		\IF{$g(x^N) \leq \varepsilon$}
		\STATE $M_N := \left\|\nabla f(x^N,\xi^{N})\right\|_*,$
		\STATE $h_N = \Theta_0 \left(\sum\limits_{t = 0}^N M_t^2\right)^{-1/2},$
		\STATE $x^{N+1}:=Mirr_{x^N}\left(h_N\nabla f(x^N,\xi^{N})\right),$  "productive step"
		\STATE$N\rightarrow I$
		\ELSE
		\STATE  (i.e. $g_{j(N)}(x^N) > \varepsilon$ \, for some $j(N) \in \{1,\ldots,m \}$) 
		\STATE $M_N := \left\|\nabla g_{j(N)}(x^N,\zeta^N)\right\|_*,$
		\STATE $h_N = \Theta_0 \left(\sum\limits_{t = 0}^N M_t^2\right)^{-1/2},$
		\STATE $x^{N+1}:=Mirr_{x^N}\left(h_N\nabla g_{j(N)}(x^N,\zeta^N)\right),$    "non-productive step"
		\ENDIF
		\STATE $N\leftarrow N+1$
		\UNTIL{$N  \geq \frac{2 \Theta_0}{\varepsilon}\left(\sum\limits_{t = 0}^{N-1} M_t^2\right)^{1/2}.$}
		\ENSURE $\bar{x}^N:=\frac{1}{N_I}\sum\limits_{k\in I}x^k.$
	\end{algorithmic}
\end{algorithm}

Denote
\begin{equation*}
\delta_k =
\begin{cases}
\langle \nabla f(x^k, \xi^k) - \nabla f(x^k), x^k-x_* \rangle \; \text{if} \; k \in I,\\
\langle \nabla g(x^k, \zeta^k) - \nabla g(x^k), x^k-x_* \rangle \; \text{if} \; k \in J.
\end{cases}
\end{equation*}

By Lemma \ref{lem1}, with $y= x^k, z= x^{k+1}$ and $x =x_*$, we have for all $k \in I$
\begin{equation} \label{ineq_for_prod_1}
\begin{aligned}
f(x^k)-f(x_*) \leq \frac{h_k}{2}\left\|\nabla f(x^k,\xi^{k})\right\|_*^2 +& \frac{V_{x^k}(x_*)}{h_k}-\frac{V_{x^{k+1}}(x_*)}{h_k} +\\
&+  \langle \nabla f(x^k, \xi^k) - \nabla f(x^k), x^k-x_*  \rangle,
\end{aligned}
\end{equation}
the same for all $k \in J$, we have (remember that, with $g_{j(k)}(\cdot)$ we mean any constraint, such that $g_{j(k)}(x^k) > \varepsilon$ ),
\begin{equation} \label{ineq_for_nonprod_2}
\begin{aligned}
g_{j(k)}(x^k)-g_{j(k)}(x_*) \leq \frac{h_k}{2}\|\nabla & g_{j(k)}(x^k,\zeta^k)\|_*^2  + \frac{V_{x^{k}}(x_*)}{h_k}-\frac{V_{x^{k+1}}(x_*)}{h_k} +\\
&+\langle \nabla g_{j(k)}(x^k, \zeta^k) - \nabla g_{j(k)}(x^k), x^k-x_* \rangle.
\end{aligned}
\end{equation}

Taking summation, in each side of \eqref{ineq_for_prod_1} and \eqref{ineq_for_nonprod_2}, over productive and non-productive steps, we get
\begin{equation*}
\begin{aligned}
\sum\limits_{k \in I}  \big( f(x^k) - f(x_{*}) \big) + \sum\limits_{k\in J}   \big( g_{j(k)}(x^k)& - g_{j(k)}(x_{*}) \big) \leq \sum\limits_{k = 0}^{N-1} \frac{h_k M_k^2}{2} +\\
& + \sum\limits_{k = 0}^{N-1} \frac{1}{h_k} \left( V_{x^{k}}(x_*) - V_{x^{k+1}}(x_*)  \right) + \sum\limits_{k=0}^{N-1} \delta_k.
\end{aligned}
\end{equation*}

Using \eqref{cond_on_Bregman}, we have
\begin{equation}
\begin{aligned}
\sum\limits_{k = 0}^{N-1} \frac{1}{h_k} \big( V_{x^{k}}(x_*) &- V_{x^{k+1}}(x_*)  \big) = \frac{1}{h_0} V_{x^{0}}(x_*) +\\
&+\sum\limits_{k=0}^{N-2}  \left[ \Big(\frac{1}{h_{k+1}} - \frac{1}{h_k} \Big) V_{x^{k+1}}(x_*) - \frac{1}{h_{N-1}} V_{x^{k}}(x_*) \right] \leq \\
& \leq  \frac{\Theta_0^2}{h_0} + \Theta_0^2 \sum\limits_{k=0}^{N-2} \Big(\frac{1}{h_{k+1}} - \frac{1}{h_k} \Big) = \frac{\Theta_0^2}{h_{N-1}}.
\end{aligned}
\end{equation}

Whence, by the definition of step-sizes $h_k$
\begin{equation}\label{111}
\begin{aligned}
\sum\limits_{k \in I}  \big(f(x^k) - f(x_*) \big) +& \sum\limits_{k\in J}  \big(g_{j(k)}(x^k) - g_{j(k)}(x_*) \big) \leq  \sum\limits_{k=0}^{N-1} \frac{\Theta_0}{2} \frac{M_k^2}{\left( \sum_{i=0}^k M_i^2 \right)^{1/2}} +\\ &+\Theta_0\left( \sum_{k=0}^{N-1} M_k^2 \right)^{1/2} + \sum\limits_{k=0}^{N-1} \delta_k \leq \\
&\leq 2 \Theta_0\left( \sum_{k=0}^{N-1} M_k^2 \right)^{1/2} + \sum\limits_{k=0}^{N -1} \delta_k,
\end{aligned}
\end{equation}
where we used the inequality
$$
	\sum\limits_{k=0}^{N-1}  \frac{M_k^2}{\left(\sum_{i=0}^k M_i^2 \right)^{1/2}} \leq 2\left( \sum_{k=0}^{N-1} M_k^2 \right)^{1/2},
$$
which can be proved by induction. Since, for $k\in J$, $g_{j(k)}(x^k) - g_{j(k)}(x_{*}) \geq g_{j(k)}(x^k) > \varepsilon$, we get
$$
	\sum\limits_{k\in J}  \big(g_{j(k)}(x^k) - g_{j(k)}(x_*) \big) > \sum\limits_{k\in J} \varepsilon = \varepsilon N_J.
$$

Thus from \eqref{111} and the stopping criterion of Algorithm \ref{algorithm2}, we have
\begin{equation}\label{1111}
\begin{aligned}
\sum\limits_{k \in I}  \big(f(x^k) - f(x_*) \big) & < 2 \Theta_0\left( \sum_{k=0}^{N-1} M_k^2 \right)^{1/2}  + \sum\limits_{k=0}^{N -1} \delta_k - \varepsilon N_J \\ & \leq \varepsilon(N_I + N_J) - \varepsilon N_J + \sum\limits_{k=0}^{N -1} \delta_k \\ & = \varepsilon N_I + \sum\limits_{k=0}^{N -1} \delta_k. 
\end{aligned}
\end{equation}

We can rewrite \eqref{1111} as follows
\begin{equation}\label{11111}
\begin{aligned}
\sum\limits_{k \in I}  f(x^k) - N_I f(x_*)   < \varepsilon N_I + \sum\limits_{k=0}^{N -1} \delta_k. 
\end{aligned}
\end{equation}

By the convexity of $f$, we get
\begin{equation}
N_I \left[ f\left( \frac{1}{N_I} \sum\limits_{k \in I} x^k \right)  - f^*  \right] < \varepsilon N_I + \sum\limits_{k=0}^{N -1} \delta_k,
\end{equation}
where $f^* = f(x_*)$. By the definition of $\bar{x}^N$ (see the Ensure of Algorithm \ref{algorithm2}), we get the following inequality
\begin{equation}\label{123}
	N_I \left( f\big( \bar{x}^N \big)  - f^*  \right) < \varepsilon N_I + \sum\limits_{k=0}^{N -1} \delta_k.
\end{equation}

As long as the inequality \eqref{123} is strict, the case of $I= \emptyset$ is impossible (i.e. $N_I \neq 0$). Now by taking the expectation in \eqref{123}  we obtain
\begin{equation*}
\mathbb{E}\left[f\big( \bar{x}^N \big)\right] - f(x_*) \leq \varepsilon + \sum\limits_{k=0}^{N-1}\mathbb{E}\left[\frac{\delta_k}{N_I}\right],
\end{equation*}
but $\sum\limits_{k=0}^{N-1}\mathbb{E}\left[\frac{\delta_k}{N_I}\right] = 0$, (see \cite{bayandina_stoc}). Thus
\begin{equation}\label{which_want_1}
\mathbb{E}\left[f\big( \bar{x}^N \big)\right] - f(x_*)\leq \varepsilon.
\end{equation}

At the same time, for $k \in I$ it holds that $g(x^k) \leq \varepsilon$. Then, by the definition of $\bar{x}^N$ and the convexity of $g$ we get
$$
	g\big(\bar{x}^N\big) \leq \frac{1}{N_I} \sum\limits_{k \in I} g(x^k) \leq \varepsilon.
$$

Thus we have come the following result

\begin{theorem}\label{theo_modification}
	Let equalities \eqref{bound_stoc_f_g_1} and inequalities \eqref{bound_stoc_f_g_2} hold. Assume that a known constant $\Theta_0>0$ is such that inequality \eqref{cond_on_Bregman} holds. Then Algorithm \ref{algorithm2} stops after no more than\\
	\begin{equation}\label{estimate_alg2}
	N = \left\lceil \frac{4 \max\{M_f^2, M_g^2\} \Theta_0^2}{\varepsilon^2}\right\rceil
	\end{equation}
	iterations and $\bar{x}^N$ is an expected $\varepsilon$-solution to problem \eqref{problem} in the sense of \eqref{expected_sol}.
\end{theorem}

\begin{remark}
From the estimate \eqref{estimate_alg2}	 we can see that Algorithm \ref{algorithm2}  achieves the complexity of the order $O \left(\varepsilon^{-2}\right)$, which is an optimal, for the studied class of non-smooth functions, from the point of view of the theory of lower bounds of estimates, according to Nemirovski and Yudin (see \cite{nemirovsky1983problem}).
\end{remark}

\section{Numerical Experiments}\label{section_numerical}
In order to compare Algorithms \ref{algorithm1} and \ref{algorithm2}, and to show the advantages  of the proposed modified algorithm some numerical tests were carried out. We consider some different examples of the following non-smooth finite-sum problem 
\begin{equation}\label{finite_sum_problem}
\min_{x} \left\{ f(x) := \frac{1}{N}\displaystyle \sum_{i=1}^{N} f_i(x) \right\},
\end{equation}
where each summand $f_i$ is a Lipschitz-continuous function. This problem is ubiquitous in many areas and applications, in particular in machine learning applications, $f$ is the total loss function whereas each $f_i$ represents the loss due to the $i$-th training sample   \cite{paper:Ding_lower_bound_2019,paper:MISO_Richtarik}. 

In our experiments, we consider the following two examples of the problem \eqref{finite_sum_problem}
\begin{example}\label{sum_linears}
	$$
		f(x) = \frac{1}{N} \displaystyle \sum_{i=1}^{N} \left| \langle a_i , x \rangle - b_i \right|, 
	$$
	where the coefficients $a_i \in \mathbb{R}^{n}$ and $b_i \in \mathbb{R}$  for each $ i = 1,\ldots,N$.

\end{example}

\begin{example}\label{sum_quad}
	$$
		f(x) = \frac{1}{N} \displaystyle \sum_{i=1}^{N} 0.5  \langle C_i x , x \rangle,
	$$  
	where  $C_i \in \mathbb{R}^{n \times n}$, for each $i = 1,\ldots,N $,  are positive definite matrices, i.e. $C_i \succ 0$.
\end{example}

For the coefficients $a_i \in \mathbb{R}^n$ and constants $b_i \in \mathbb{R}$ $(i=1,\ldots,N)$, in example \ref{sum_linears}, with different values of $N$. Let $A \in \mathbb{R}^{N \times (n+1)} $ be a matrix with entries drawn from different random distributions. Then $a_i^T$ are rows in the matrix $A' \in \mathbb{R}^{N \times n} $, which is obtained from $A$, by eliminating the last column, and $b_i$ are the entries of the last column in the matrix $A$.   The positive definite matrices $ C_i \succ 0$ $(i=1,\ldots,N)$, in example \ref{sum_quad}, with different values of $N$, are drawn from different random distributions. In more details, the entries of $A$ and $C_i( i=1,\ldots,N)$, with different values of $N$, are drawn
\begin{enumerate}
	\item When $N = 75$, from the Gumbel distribution with the location of the mode equaling $1$ and the scale parameter equaling $2$.
	
	\item When $N = 100$, from the standard exponential distribution with a scale parameter of $1$.
	
	\item When $N = 150$, from the uniform distribution over $[0,1)$.
\end{enumerate}

For the functional  constraint $g(x) = \max\limits_{i \in \overline{1,  m}} \{g_i(x)\}$, we take $m = 50, n=1500$ and  $g_i(x) =\langle \alpha_i , x \rangle + \beta_i$ linear functionals, where the coefficients $\alpha_i \in \mathbb{R}^n$ and $\beta_i \in \mathbb{R}$ for $i=1,\ldots,m$  are taken as follows: Let $B \in \mathbb{R}^{m \times (n+1)} $ be a Toeplitz matrix  with the first row $(1,1,\ldots,1) \in \mathbb{R}^{n+1}$ and the first column $(1,2,\ldots,m)^T$. Then $\alpha_i^T$ are rows in the matrix $B' \in \mathbb{R}^{m \times n} $, which is obtained from $B$, by eliminating the last column, and $\beta_i$ are the entries of the last column in the matrix $B$, i.e. the eliminated column. 

For more clarification, when $m=10$ and  $n=14$, then the Toeplitz matrix $B$ with the first row $(1,1,\ldots,1) \in \mathbb{R}^{15}$ and the first column $(1,2,\ldots,10)^T$ has the form \\
$$ B =
\left(
\begin{array}{ccccccccccccccc}
1\, & 1\, & 1\, & 1\, &1\, & 1\, & 1\, & 1\, & 1\, & 1\,&1\,& 1\, & 1\, & 1\,&1\\
2 &  1 & 1 & 1  & 1  & 1  & 1  & 1  & 1  & 1 &1& 1  & 1  & 1 &1\\
3 & 2 & 1 & 1 & 1 & 1 & 1 & 1 & 1 & 1 &1& 1  & 1  & 1 &1\\
4 & 3 & 2 & 1 & 1 & 1 & 1 & 1 & 1 & 1 &1& 1  & 1  & 1 &1\\
5 & 4 & 3 & 2 & 1 & 1 & 1 & 1 & 1 & 1&1& 1  & 1  & 1 &1\\
6 & 5 & 4 & 3 & 2 & 1 & 1 & 1 & 1 & 1 &1& 1  & 1  & 1 &1\\
7 & 6 & 5 & 4 & 3 & 2 & 1 & 1 & 1 & 1 &1& 1  & 1  & 1 &1\\
8 & 7 & 6 & 5 & 4 & 3 & 2 & 1 & 1 & 1&1& 1  & 1  & 1 &1\\
9 & 8 & 7 & 6 & 5 & 4 & 3 & 2 & 1 & 1 &1& 1  & 1  & 1 &1\\
10 & 9 & 8 & 7 & 6 & 5 & 4 & 3 & 2  & 1 &1& 1  & 1  & 1 &1
\end{array}
\right) \in \mathbb{R}^{10 \times 15}.
$$

The proximal structure is given by Euclidean norm and squared Euclidean norm as a prox-function. We choose starting point $x^0 = \left(\frac{1}{\sqrt{n}}, \frac{1}{\sqrt{n}},\ldots,\frac{1}{\sqrt{n}}\right)\in \mathbb{R}^{n}$, $\varepsilon = 0.05$, and $Q = \{x = (x_1, x_2, \ldots , x_{n}) \in \mathbb{R}^{n}\,|\, x_1^2 + x_2^2 + \ldots + x_{n}^2 \leq 1\}.$

For any $x = (x_1,\ldots , x_{n})$ and $y = (y_1, \ldots , y_{n})$ in $Q$, the following inequality holds
$$
	\frac{1}{2}\|x-y\|_2^{2} = \frac{1}{2}\sum\limits_{k=1}^{n}(x_k - y_k)^{2} \leq x_1^2 + \ldots + x_{n}^2 + y_1^2 + \ldots + y_{n}^2 \leq 2.
$$
Therefore, we can choose $\Theta_0 = \sqrt{2}.$

Our experiments are motivated by the need to solve the problem \eqref{problem} when either the dimension $n$ is large or when the objective function $f$ is of a finite sum structure, as in examples \ref{sum_linears} and \ref{sum_quad}, with  $N$, the number of components, being large.

We run Algorithms \ref{algorithm1} and \ref{algorithm2}, in order to both Examples \ref{sum_linears} and \ref{sum_quad}, with $m = 50, n=1500$. The results of the work of Algorithms \ref{algorithm1} and \ref{algorithm2} are represented in Table \ref{tab:1_2}, below. These results demonstrate the comparison between the number of iterations and the running time (in seconds) for each algorithm.

All experiments were implemented in Python 3.4, on a computer fitted with Intel(R) Core(TM) i7-8550U CPU @ 1.80GHz, 1992 Mhz, 4 Core(s), 8 Logical Processor(s). RAM of the computer is 8GB.

\begin{table}[]
	\centering
	\caption{Results of Algorithms \ref{algorithm1} and \ref{algorithm2}, for Examples \ref{sum_linears} and  \ref{sum_quad}, in $\mathbb{R}^{1500}$.}
	\label{tab:1_2}
	\begin{tabular}{|c|c|c|c|c|}
		\hline
		\multirow{3}{*}{} & \multicolumn{4}{c|}{Example \ref{sum_linears}}                         \\ \cline{2-5} 
	$N$	& \multicolumn{2}{c|}{Algorithm \ref{algorithm1}} & \multicolumn{2}{c|}{Algorithm \ref{algorithm2}} \\ \cline{2-5} 
		& Iterations & Time (sec)  & Iterations  & Time (sec)  \\ \hline
		$75$ &30\,157 &618.79 &27\,007 &22.47\\ 
		$100$&12\,827 &254.34 &11\,071 & 10.04 \\ 
		$150$&7\,452  &139.99 &5\,713  &4.62  \\ \hline \hline
		\multirow{3}{*}{} & \multicolumn{4}{c|}{Example  \ref{sum_quad}}                        \\ \cline{2-5} 
	$N$	& \multicolumn{2}{c|}{Algorithm \ref{algorithm1}} & \multicolumn{2}{c|}{Algorithm \ref{algorithm2}} \\ \cline{2-5} 
		& Iterations & Time (sec)  & Iterations  & Time (sec) \\ \hline
		$75$  &104\,513  &2008.12 &90\,154  &82.38  \\ 
		$100$ &18\,814   &358.02  &17\,584  &15.3 \\ 
		$150$ &5\,451    &115.47  &4\,834   &5.45  \\ \hline
	\end{tabular}
\end{table}

From Table \ref{tab:1_2}, in order to both examples \ref{sum_linears} and \ref{sum_quad}, we can see that the modified Algorithm \ref{algorithm2} always works better than Algorithm \ref{algorithm1}. It is clearly shown in all experiments according to the number of iterations and especially according to the running time of the algorithms. The running time of Algorithm \ref{algorithm2} is very small  compared to the running time of Algorithm \ref{algorithm1} (on average, it is smaller 25 times). This feature of the Algorithm \ref{algorithm2} is very important in all applications of mathematical optimization.  

\begin{remark}
Now, as in the previous, to compare  Algorithms \ref{algorithm1} and \ref{algorithm2}, with $m = 50, n = 100$ and different values of $N$, some additional numerical tests were carried out. The coefficients $\alpha_i \in \mathbb{R}^n$ and $\beta_i \in \mathbb{R}$, for each $i=1, \ldots,m$, are the entries of the Toeplitz matrix, which is described above. The entries of the matrices $C_i( i=1,\ldots,N)$ are drawn from the uniform distribution over $[0,1)$.  We run Algorithms \ref{algorithm1} and \ref{algorithm2} with the same previous parameters $\varepsilon = 0.05, \Theta_0 = \sqrt{2}$ and the set $Q$.  The results of Algorithms \ref{algorithm1} and \ref{algorithm2}, in order to the examples \ref{sum_linears} and \ref{sum_quad}  are represented in Table \ref{tab:N_1_2}, below.  These results demonstrate the comparison between the number of iterations and the running time (in seconds) for each algorithm, with different values of $N$.

From Table \ref{tab:N_1_2}, we can see that Algorithm \ref{algorithm2} works better than Algorithm \ref{algorithm1} according to the number of iterations and especially according to the running time of algorithms.

\begin{table}[t]
	\centering
	\caption{The results of Algorithms \ref{algorithm1} and \ref{algorithm2}, for Examples \ref{sum_linears} and  \ref{sum_quad}, with different values of $N$.}
	\label{tab:N_1_2}
	\begin{tabular}{|c|c|c|c|c|}
		\hline
		\multirow{3}{*}{} & \multicolumn{4}{c|}{Example \ref{sum_linears}}                         \\ \cline{2-5} 
		$N$	& \multicolumn{2}{c|}{Algorithm \ref{algorithm1}} & \multicolumn{2}{c|}{Algorithm \ref{algorithm2}} \\ \cline{2-5} 
		& Iterations & Time (sec)  & Iterations  & Time (sec)  \\ \hline
		$1\,000 $  &6\,717  &9.770  &5\,366  &0.476  \\ 
		$5\,000 $  &5\,726  &7.975  &5\,334  &0.452  \\ 
		$10\,000 $ &8\,017  &11.076 &5\,574  &0.500  \\ 
		$15\,000 $ &6\,427  &8.890  &5\,243  &0.445  \\ 
		$25\,000 $ &6\,775  &9.530  &5\,348  &0.474 \\ 
		$50\,000 $ &7\,339  &10.232 &6\,187  &0.582  \\
		$75\,000 $ &6\,599  &9.160  &5\,287  &0.452  \\ 
		$100\,000$ &6\,235  &8.665 &5\,400  &0.456 \\ 
		$125\,000$ &6\,709  &9.175  &6\,095  &0.512  \\
		$150\,000$ &6\,928  &9.671  &5\,360  &0.471 \\ \hline \hline 
		\multirow{3}{*}{} & \multicolumn{4}{c|}{Example  \ref{sum_quad}}                        \\ \cline{2-5} 
		$N$	& \multicolumn{2}{c|}{Algorithm \ref{algorithm1}} & \multicolumn{2}{c|}{Algorithm \ref{algorithm2}} \\ \cline{2-5} 
		& Iterations & Time (sec)  & Iterations  & Time (sec) \\ \hline
		$1\,000 $  &6\,519  &10.496  &5\,178  &0.656  \\ 
		$2\,500 $  &6\,238  &9.750   &4\,634  &0.523  \\ 
		$5\,000 $  &5\,364  &8.287   &4\,615  &0.679  \\ 
		$7\,500 $  &5\,862  &9.255   &5\,029  &0.677  \\ 
		$10\,000 $ &6\,025  &9.331   &4\,506  &0.569  \\ 
		$12\,500 $ &5\,341  &10.687  &4\,688  &0.672 \\ 
		$15\,000 $ &6\,227  &12.981  &4\,995  &0.576  \\ 
		$17\,500 $ &5\,847  &9.509   &4\,616  &0.603  \\ 
		$20\,000 $ &5\,486  &8.515   &4\,760  &0.620 \\ 
		$22\,500 $ &6\,294  &10.140  &4\,551  &0.585  \\ 
		$25\,000 $ &6\,055  &11.598  &4\,534  &0.596 \\ \hline
	\end{tabular}
\end{table}

\end{remark}

\subsection{Additional Experiments: Fermat-Torricelli-Steiner problem}

In this subsection some additional numerical experiments connected with the analogue of the well-known Fermat-Torricelli-Steiner problem with some non-smooth functional constraints, were carried out. 

For a given  set  $\{A_k=(a_{1k},a_{2k},\ldots,a_{nk}); \, k=\overline{1,N}\}$ of $N$ points, in $n$-dimensional Euclidean space $\mathbb{R}^n$, we need to solve the problem \eqref{problem}, where the objective function $f$ is given by 
$$
	f(x):=\sum\limits_{k=1}^N\sqrt{(x_1-a_{1k})^2+(x_2-a_{2k})^2+\ldots+(x_n-a_{nk})^2} = \sum\limits_{k=1}^N \|x - A_k\|_2.
$$

The functional constraint is given by $g(x) = \max\limits_{i \in \overline{1,  m}} \{g_i(x) =\langle \alpha_i , x \rangle + \beta_i\} $, where the coefficients $\alpha_i \in \mathbb{R}^n$ and $\beta_i \in \mathbb{R}$ are taken as in the previous experiments (the entries of the Toeplitz matrix $B$).

We take the points $A_k (k \in \overline{1,N})$  in the unit ball $Q$. The coordinates of these points are drawn from the uniform distribution over $[0,1)$.

We choose the standard Euclidean proximal setup,  starting point $x^0 = \textbf{0} \in \mathbb{R}^n$ and $\Theta_0 = \sqrt{2}$. We run Algorithms \ref{algorithm1} and \ref{algorithm2} with  $n = 1000, m = 250, N = 100$ and different values of accuracy $\varepsilon \in \{ 1/2^i: i=1,2,3,4,5,6\}$.

\begin{figure}[htb]
	\centering
	\begin{subfigure}{.5\textwidth}
		\centering
		\includegraphics[width=1\linewidth]{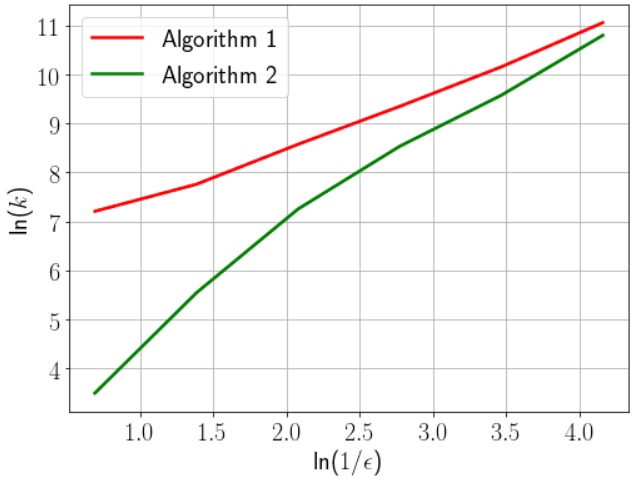}
		\caption{}
		\label{results:Fermat_iter}
	\end{subfigure}%
	\begin{subfigure}{.5\textwidth}
		\centering
		\includegraphics[width=1\linewidth]{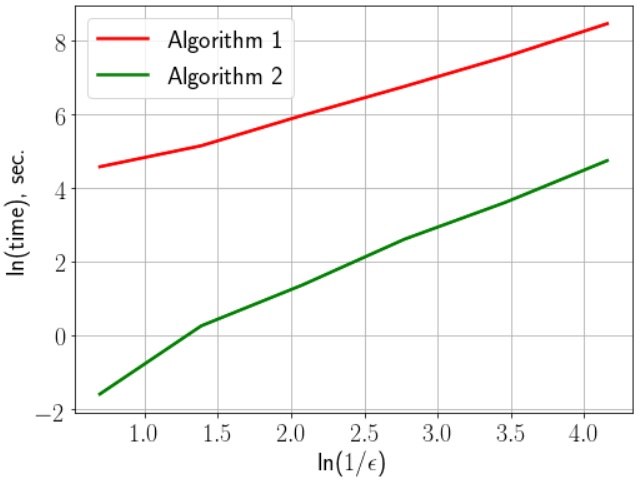}
		\caption{}
		\label{results:Fermat_time}
	\end{subfigure}
	\caption{The results of Algorithms \ref{algorithm1} and \ref{algorithm2}, for Fermat-Torricelli-Steiner problem.}
	\label{fig_exp10_m1}
\end{figure}

The results of the work of  Algorithms \ref{algorithm1} and \ref{algorithm2}, are presented in Fig. \ref{results:Fermat_iter} (the number of iterations produced by the studied algorithms to reach an $\varepsilon$-solution of the proposed problem as a function of accuracy) and Fig. \ref{results:Fermat_time} (the required running time of the studied algorithms, in seconds, as a function of accuracy).

From Fig. \ref{results:Fermat_iter} and  Fig. \ref{results:Fermat_time}, we see that both Algorithms \ref{algorithm1} and \ref{algorithm2} are optimal, where they achieve the complexity of the order $O \left(\varepsilon^{-2}\right)$, which is optimal estimate for the studied  class of non-smooth functions. But Algorithm \ref{algorithm2} is more efficiently and works better than Algorithm \ref{algorithm1}, according to the number of iterations and the running time. We note that the running time of Algorithm \ref{algorithm1} is very long compared with the running time of Algorithm \ref{algorithm2}, where by Algorithm \ref{algorithm2} one needs a few seconds, when needs more and more minutes by Algorithm \ref{algorithm1}, to achieve a solution and to reach its stopping criterion. Therefore the efficiency of Algorithm \ref{algorithm2} is represented by its very high execution speed compared with Algorithm \ref{algorithm1}.

\section{Conclusions}
In this work, a new modification of an adaptive stochastic mirror descent algorithm was proposed to solve the stochastic setting of the convex minimization problem in the case of Lipschitz-continuous objective function and several convex functional constraints. In each iteration of the proposed modified algorithm, we calculate the stochastic (sub)gradient of the objective function or the functional of constraint, which is prevalent and effective in Machine Learning scenarios, large-scale optimization problems, and their applications. The proposed modification allows saving the running time of algorithm  due to the consideration of not all functional constraints on non-productive steps. Furthermore, it has been proved a theorem to estimate the rate of convergence of the proposed modified algorithm. Numerical experiments for a geometrical problem, Fermat-Torricelli-Steiner problem, with convex constraints are presented. The results of carried out numerical experiments illustrate the advantages of the modified Algorithm \ref{algorithm2} and illustrate that the running time of this Algorithm is very small compared to the running time of the standard Algorithm \ref{algorithm1}.

\textbf{Acknowledgments:}
The author is very grateful to Alexander V. Gasnikov, Fedor S. Stonyakin and Alexander G. Biryukov for  fruitful discussions.

\end{document}